\documentstyle[12pt]{article}
\catcode`\@=11
\@addtoreset{equation}{section}

\catcode`\@=12
\newtheorem{Theorem}{Theorem}[section]
\newtheorem{Definition}{Definition}[section]
\newtheorem{Proposition}{Proposition}[section]
\newtheorem{Lemma}{Lemma}[section]
\newtheorem{Corollary}{Corollary}[section]
\newtheorem{Note}{Note}[section]
\title{$J-$holomorphic Curves, Legendre Submanifolds and Reeb Chords
\thanks{Project 19871044 Supported by NSF}}

\author{Renyi Ma \\
Department of Mathematics \\
Tsinghua University \\
Beijing, 100084\\
People's Republic of China\\
rma@math.tsinghua.edu.cn}

\date { }

\begin{document}
\textwidth=165mm
\textheight=185mm
\parindent=8mm
\frenchspacing
\maketitle

\begin{abstract}
In this article, we prove that there
exists at least one chord which is characteristic
of Reeb vector field connecting a given Legendre submanifold
in a closed contact manifold with any contact form.
\end{abstract}

\noindent{\bf 2000 MR Subject Classification} 32Q65, 53D35,53D12

\section{Introduction and results}

Let $\Sigma$ be a smooth closed oriented manifold of dimension
$2n-1$. A contact form on $\Sigma$ is a $1-$form such that
$\lambda \wedge (d\lambda )^{n-1}$ is a volume form on $\Sigma$.
Associated to $\lambda$ there are two important structures.
First of all the so-called Reed vectorfield $\dot x=X$ defined
by
$$i_X\lambda   \equiv 1, \ \ i_Xd\lambda  \equiv 0;$$
and secondly the contact structure $\xi =\xi _{\lambda }
\mapsto \Sigma $ given by
$$\xi _{\lambda }=\ker (\lambda )\subset T \Sigma .$$
By a result of Gray, \cite{gra}, the contact structure is very
stable. In fact, if $(\lambda  _t  )_{t\in  [0,1]}$ is a smooth
arc of contact forms inducing the arc of contact structures
$(\xi _t)_{t\in [0,1]}$, there exists a smooth arc
$(\psi _t)_{t\in [0,1]}$
of diffeomorphisms with $\psi _0=Id$, such that
\begin{equation}
T\Psi _t(\xi _0)=\xi _t  \label{eq:1.1}
\end{equation}
here it is importent that $\Sigma $ is compact. From (\ref{eq:1.1}) and
the fact that $\Psi _0=Id$ it follows immediately that there
exists a smooth family of maps $[0,1]\times \Sigma \mapsto
(0,\infty ):(t, m)\to f_t(m)$ such that
\begin{equation}
\Psi ^*_t\lambda _t=f_t\lambda _0
\end{equation}
In contrast to the contact structure the dynamics of the
Reeb vectorfield changes drastically under small perturbation
and in general the flows associated to $X_t$ and $X_s$ for
$t\neq s$ will not be conjugated.

Concerning the dynamics of Reeb flow,
there is a well-known conjecture raised by Arnold in \cite{ar}
which concerned the Reeb orbit and Legendre submanifold in
a contact manifold. If $(\Sigma ,\lambda )$ is a contact manifold
with contact form $\lambda $ of dimension $2n-1$, then a Legendre
submanifold is a submanifold ${\cal L}$ of $\Sigma $, which is
$(n-1)$dimensional and everywhere tangent to the contact
structure $\ker \lambda $. Then a characteristic
chord for $(\lambda ,{{\cal {L}}})$ is a smooth path
$$x:[0,T]\to M,T>0$$
with
$$\dot x(t)=X_{\lambda }(x(t)) \ for \ t\in(0,T),$$
$$x(0),x(T)\in {\cal {L}}$$
Arnold raised the following conjecture:

\vskip 3pt

{\bf Conjecture}(see\cite{ar}). Let $\lambda _0$ be the standard tight
contact form
$$\lambda _0={{1}\over {2}}(x_1dy_1-y_1dx_1+x_2dy_2-y_2dx_2)$$
on the three sphere
$$S^3=\{ (x_1,y_1,x_2,y_2)\in R^4|x_1^2+y_1^2+x_2^2+y_2^2=1\}.$$
If $f:S^3\to (0,\infty )$ is a smooth function
and ${\cal {L}}$ is a Legendre
knot in $S^3$, then
there is a characteristic chord for $(f\lambda _0,{\cal {L}})$.

\vskip 3pt

The main results of this paper as following:

\begin{Theorem}
Let $(\Sigma ,\lambda )$ be a contact manifold with
contact form $\lambda $, $X_{\lambda } $ its Reeb vector
field, ${\cal {L}}$ a closed Legendre submanifold, then
there exists at least one characteristic chord for
$(X_\lambda ,{\cal {L}})$.
\end{Theorem}

\begin{Corollary}
(\cite{ma,mo}) Let $(S^3,f\lambda _0)$ be a tight contact manifold
with contact form $f\lambda _0$, $X_{f\lambda _0} $ its Reeb
vector field, ${\cal {L}}$ a closed Legendre submanifold, then
there exists at least one characteristic chord for $(X_{f\lambda
_0},{\cal {L}})$.
\end{Corollary}

{\bf Sketch of proofs}: We work in the framework
as in \cite{al,gro, ma}. In Section 2, we study the
linear Cauchy-Riemann operator and sketch some basic
properties. In section 3, first we construct a
Lagrangian submanifold $W$ under the assumption that
there does not exists Reeb chord conneting the Legendre
submanifold $\cal {L}$; second, we study
the space ${\cal D}(V,W)$
of contractible disks in manifold $V$ with boundary
in Lagrangian submanifold $W$ and construct a Fredholm
section of tangent bundle of  ${\cal  D }(V,W)$.
In section 4, following \cite{al,gro, ma}, we prove that the Fredholm
section is not proper by using an
special anti-holomorphic section as in \cite{al,gro, ma}.  In section 5-6,
we use a geometric argument
to prove the boundaries of $J-$holomorphic
curves remain in a finite part of Lagrangian
submanifold $W$.
In the final section, we
use nonlinear Fredholm trick in \cite{al,gro, ma}
to
complete our proof.

\section{Linear Fredholm Theory}

For $100<k<\infty $ consider the Hilbert space
$V_k$ consisting of all maps $u\in H^{k,2}(D, C\times C^n)$,
such that $u(z)\in \{izR\}\times R^n\subset C\times C^n$
for almost all $z\in
\partial D$. $L_{k-1}$ denotes the usual Sobolev
space
$H_{k-1}(D, C\times C^n)$. We define an operator
$\bar \partial :V_k\mapsto L_{k-1}$ by
\begin{equation}
\bar \partial u=u_s+iu_t
\end{equation}
where the coordinates on $D$ are
$(s,t)=s+it$, $D=\{ z||z|\leq 1\}  $.
The following result is well known(see\cite{wen}).
\begin{Proposition}
$\bar \partial :V_k\mapsto L_{k-1}$ is a surjective real linear
Fredholm operator of index $n+3$. The kernel consists of
$(a_0+isz-\bar a_0z^2,s_1,...,s_n)$, $a_0\in C$, $s,s_1,...,s_n\in
R$.
\end{Proposition}
Let $(C^n, \sigma =-Im(\cdot ,\cdot ))$ be the standard
symplectic space. We consider a real $n-$dimensional plane
$R^n\subset C^n$. It is called Lagrangian if the
skew-scalar product of any two vectors of $R^n$ equals zero.
For example, the plane $\{ (p,q)|p=0\}$ and $\{ (p,q)|q=0\}$ are
two transversal Lagrangian
subspaces. The manifold of all (nonoriented) Lagrangian subspaces of
$R^{2n}$ is called the Lagrangian-Grassmanian $\Lambda (n)$.
One can prove that the fundamental group of
$\Lambda (n)$ is free cyclic, i.e.
$\pi _1(\Lambda (n))=Z$. Next assume
$(\Gamma (z))_{z\in \partial D}$ is a smooth map
associating to a point $z\in \partial D$ a Lagrangian
subspace $\Gamma (z)$ of $C^n$, i.e.
$(\Gamma (z))_{z\in \partial D}$ defines a smooth curve
$\alpha $ in the Lagrangian-Grassmanian manifold $\Lambda (n)$.
Since $\pi _1(\Lambda (n))=Z$, one have
$[\alpha ]=ke$, we call integer $k$ the Maslov index
of curve $\alpha $ and denote it by $m(\Gamma )$, see(\cite{ag,wen}).

Now let $z:S^1\mapsto \{R\times R^n\subset C\times C^n\}\in \Lambda (n+1)$ be a constant curve.
Then
it defines a constant loop $\alpha $ in Lagrangian-Grassmanian
manifold $\Lambda (n+1)$. This loop defines
the Maslov index $m(\alpha )$ of the map
$z$ which is easily seen to be zero.

  Now Let $(V,\omega )$ be a symplectic manifold,
$W\subset V$ a closed Lagrangian submanifold. Let
$(\bar V,\bar \omega )=(D\times V, \omega _0+\omega )$ and
$\bar W=\partial D\times W$. Let
$\bar u=(id,u):(D,\partial D)\to (D\times V,\partial D\times W)$
be a smooth map homotopic to the map
$\bar u_0=(id,u_0)$, here $u_0:(D, \partial D)\to p\in W\subset V$.
Then $\bar u^*TV$ is a symplectic vector bundle on $D$ and
$(\bar u|_{\partial D})^*T\bar W$ be a Lagrangian subbundle in
$\bar u^*T\bar V|_{\partial D}$. Since $\bar u:(D,\partial D)\to
(\bar V,\bar W)$ is homotopic to
$\bar u_0$, i.e., there exists
a homotopy $h:[0,1]\times (D,\partial D)\to (\bar V,\bar W)$ such
that $h(0,z)=(z,p),h(1,z)=\bar u(z)$, we can take
a trivialization of the symplectic
vector bundle $h^*T\bar V$ on $[0,1]\times (D,\partial D)$ as
$$\Phi (h^*T\bar V)=[0,1]\times D\times C\times C^n$$
and
$$\Phi ((h|_{[0,1]\times \partial D})^*T\bar W)\subset [0,1]\times
S^1\times C\times C^n$$
Let
$$\pi _2: [0,1]\times D\times C\times C^n\to C\times C^n$$
then
$$\tilde h: (s,z)\in [0,1]\times S^1
\to \pi _2\Phi (h|_{[0,1]\times \partial D})^*T\bar W|(s,z)\in \Lambda (n+1).$$

\begin{Lemma}
Let $\bar u: (D,\partial D) \rightarrow (\bar V,\bar W)$
be a $C^k-$map $(k\geq 1)$
as above. Then,
$$m(\tilde u)=2.$$
\end{Lemma}
Proof.  Since $\bar u$ is homotopic to
$\bar u_0$ in $\bar V$ relative to $\bar W$,
by the above argument we have
a homotopy $\Phi _s$ of trivializations
such that
$$\Phi _s(\bar u^*TV)=D\times C\times C^n$$
and
$$\Phi _s((\bar u|_{\partial D})^*T\bar W)\subset S^1\times C\times C^n$$
Moreover
$$\Phi _0(\bar u|_{\partial D})^*T\bar W=S^1\times {izR}\times R^n$$
So, the homotopy induces
a homotopy $\tilde h$ in Lagrangian-Grassmanian
manifold. Note that $m(\tilde h(0, \cdot ))=0$.
By the homotopy invariance of Maslov index,
we know that $m(\tilde  u|_{\partial D})=2$.

\vskip 5pt

   Consider the partial differential equation
\begin{eqnarray}
&&\bar \partial \bar u+A(z)\bar u=0  \ on \ D  \cr
&&\bar u(z)\in \Gamma (z) ({izR}\times R^n)\ for \ z\in \partial D \cr
&&\Gamma (z)\in GL(2(n+1),R)\cap Sp(2(n+1))\cr
&&m(\Gamma )=2 \ \ \ \ \ \ \ \
\end{eqnarray}

For $100<k<\infty $ consider the Banach space $\bar V_k $
consisting of all maps $u\in H^{k,2}(D, C^n)$ such
that  $u(z)\in \Gamma (z)$ for almost all $z\in
\partial D$. Let $L_{k-1}$ the usual Sobolev space $H_{k-1}(D,C\times C^n)$

\begin{Proposition}
$\bar \partial : \bar V_k \rightarrow L_{k-1}$
is a real linear Fredholm operator of index n+3.
\end{Proposition}

\section{Nonlinear Fredholm Theory}

\subsection{Constructions of Lagrangian submanifolds}

Let $(\Sigma ,\lambda )$ be a contact manifolds with contact form
$\lambda $ and $X$ its Reeb vector field, then
$X$ integrates to a Reeb flow $\eta _t$ for $t\in R^1$.
Consider the form $d(e^a\lambda )$
on the manifold
$(R\times \Sigma )$, then one can check
that $d(e^a\lambda )$ is a symplectic
form on $R\times \Sigma $. Moreover
One can check that
\begin{eqnarray}
&&i_X(e^a\lambda )=e^a \\
&&i_X(d(e^a\lambda ))=-de^a
\end{eqnarray}
So, the symplectization of Reeb vector field $X$ is the Hamilton
vector field of $e^a$ with respect to the symplectic form
$d(e^a\lambda )$. Therefore the Reeb flow lifts to the Hamilton
flow $h_s$ on $R\times \Sigma $(see\cite{ag}).

Let ${\cal L}$ be a closed Legendre submanifold
in $(\Sigma ,\lambda )$, i.e.,
there exists a smooth embedding $Q:{\cal L}\to \Sigma $ such that
$Q^*\lambda |_{\cal L}=0$, $\lambda |Q(L)=0$. We also
write ${\cal {L}}=Q({\cal {L}})$. Let
$$(V',\omega ')=(R\times \Sigma ,d(e^a\lambda ))$$
and
\begin{eqnarray}
&&W'={\cal {L}}\times R,
\ \ W'_s={\cal L}\times \{ s\} ;\cr
&&L'=(0,\cup _s \eta _s(Q({\cal {L}}))),
\ \ L'_s=(0, \eta _s(Q({\cal {L}})))
\label{eq:3.w}
\end{eqnarray}
define
\begin{eqnarray}
&&G':W'\to V'  \cr
&&G'(w')=G'(l,s)=(0,\eta _s(Q(l))) \label{eq:3.ww}
\end{eqnarray}
\begin{Lemma}
There does not exist any Reeb chord connecting Legendre
submanifold ${\cal {L}}$
in $(\Sigma ,\lambda )$ if and only if
$G'(W'_s)\cap G'(W'_{s'})$ is empty for $s\ne s'$.
\end{Lemma}
Proof. Obvious.
\begin{Lemma}
If there does not exist any Reeb chord for $(X_\lambda ,{\cal {L}})$
in $(\Sigma ,\lambda )$ then
there exists a smooth embedding
$G':W'\to V'$ with $G'(l,s)=(0,\eta _s(Q(l)))$
such that
\begin{equation}
G'_K:{\cal L}\times (-K, K)\to V'  \label{eq:3.www}
\end{equation}
is a regular open Lagrangian embedding for any finite positive $K$.
We denote $W'(-K,K)=G'_K({\cal L}\times (-K,K))$
\end{Lemma}
Proof. One check
\begin{equation}
{G'}^*(d(e^a\lambda ))=\eta (\cdot ,\cdot )^*d\lambda
=(\eta _s^*d\lambda +i_Xd\lambda \wedge ds)=0
\end{equation}
This implies that ${G}'$ is a Lagrangian embedding, this proves
Lemma3.2.

\vskip 3pt

In fact the above proof checks that
\begin{equation}
{G'}^*(\lambda )=\eta (\cdot , \cdot )^*\lambda =
\eta _s^*\lambda +i_X\lambda ds=ds.
\end{equation}
i.e., $W'$ is an exact Lagrangian submanifold.

\vskip 3pt

Now we construct an isotopy of Lagrangian embeddings as follows:
\begin{eqnarray}
&&F':{\cal L}\times R\times [0,1]\to R\times \Sigma \cr
&&F'(l,s,t)=(a(s,t),G'(l,s))=(a(s,t),\eta _s(Q(l)))  \cr
&&F'_t(l,s)=F'(l,s,t)
\end{eqnarray}
\begin{Lemma}
If there does not exist any Reeb chord for $(X_\lambda ,{\cal {L}})$
in $(\Sigma ,\lambda )$ and we choose the smooth  $a(s,t)$ such that
$\int _0^sa(\tau ,t)d\tau $ and $\int _s^0a(\tau ,t)d\tau  $ exists,
then $F'$ is an exact isotopy of Lagrangian
embeddings(not regular). Moreover if $a(s,0)\ne a(s,1)$, then
$F'_0({\cal L}\times R)\cap F'_1({\cal L}\times R)=\emptyset $.
\end{Lemma}
Proof. Let $F'_t=F'(\cdot,t):{\cal L}\times R\to R\times \Sigma $.
It is obvious that $F'_t$ is an embedding. We check
\begin{eqnarray}
F'^*(d(e^a\lambda ))
&=&d(F'^*(e^a\lambda ))\cr
&=&d(e^{a(s,t)}G'^*\lambda ) \cr
&=&d(e^{a(s,t)}ds) \cr
&=&e^{a(s,t)}(a_sds+a_tdt)\wedge ds \cr
&=&e^{a(s,t)}a_tdt\wedge ds
\end{eqnarray}
which shows that $F'_t$ is a Lagrangian embedding for fixed $t$.
Moreover for fixed $t$,
\begin{eqnarray}
F_t^*(e^a\lambda )&=&e^{a(s,t)}ds \cr
&=&\left\{ \begin{array}{ll}
d(\int _{0}^se^{a(\tau ,t)}d\tau )  &\mbox{for $s\geq 0$}\cr
d(-\int _s^0e^{a(\tau ,t)}d\tau )   &\mbox{for $s\leq 0$}
\end{array}
\right.
\end{eqnarray}
which shows that $F'_t$ is an exact Lagrangian embedding, this proves
Lemma3.3.

\vskip 5pt

Now we take
$a(s,t)={{a_0\varepsilon t}\over {8}}e^{-s^2}$ which satisfies
the assumption in Lemma 3.3, then
$$F':{\cal L}\times R\times [0,1]\to R\times \Sigma $$
$$F'(l,s,t)=(a(s,t),\eta _s(G(l))).$$
Let
\begin{equation}
\psi _0(s,t)=se^{a(s,t)}a_s=-2{{a_0\varepsilon t}\over {8}}
e^{({{a_0\varepsilon t}\over {8}}e^{-s^2})-s^2}s^2 \label{eq:fu}
\end{equation}
\begin{equation}
\psi _1(s,t)=\int _{-\infty }^s\psi _0(\tau, t)d\tau \label{eq:fu1}
\end{equation}
\begin{equation}
\psi ={{\partial \psi _1}\over {\partial t}}-se^{a(s,t)}a_t \label{eq:13}
\end{equation}
and compute
\begin{eqnarray}
{F'}^*(e^a\lambda )&=&e^{a(s,t)}ds \cr
&=&d(se^{a(s,t)} )-sde^{a(s,t)}\cr
&=&d(se^{a(s,t)})-se^{a(s,t)}a_sds-se^{a(s,t)}a_tdt  \cr
&=&d(se^{a(s,t)})-d_s{\psi _1}-se^{a(s,t)}a_tdt \cr
&=&d((se^{a(s,t)})-\psi _1)
+{{\partial \psi _1}\over {\partial t}}dt
-se^{a(s,t)}a_tdt  \cr
&=&d\Psi '+{{\partial \psi _1}\over {\partial t}}dt-se^{a(s,t)}a_tdt \cr
&=&d\Psi '-\psi (s,t)dt \cr
&=&d\Psi '-\tilde {l}'  \label{eq:14}
\end{eqnarray}

   Let
$(V',\omega ')=(R\times \Sigma ,d(e^a\lambda ))$,
$W'={\cal L}\times R$, and
$(V,\omega )=(V'\times C,\omega '\oplus \omega _0)$.
As in \cite{gro}, we use figure eight trick invented by Gromov to
construct a Lagrangian submanifold in $V$ through the
Lagrange isotopy $F'$ in $V'$.
Fix a positive $\delta <1$ and take a $C^{\infty }$-map $\rho :S^1\to
[0,1]$, where the circle $S^1$
via parametrized by $\Theta \in [-1,1]$,
such that the $\delta -$neighborhood $I_0$ of $0\in S^1$ goes to
$0\in [0,1]$ and $\delta -$neighbourhood $I_1$ of $\pm 1\in S^1$
goes $1\in [0,1]$. Let
\begin{eqnarray}
\tilde {l}&=&-\psi (s,\rho (\Theta ))\rho '(\Theta )d\Theta \cr
&=&-\Phi d\Theta
\end{eqnarray}
be the pull-back of the form
$\tilde {l}'=-\psi (s,t)dt $ to $W'\times S^1$ under the map
$(w',\Theta )\to (w',\rho (\Theta ))$ and
assume without loss of generality $\Phi $
vanishes on $W'\times (I_0\cup I_1)$.

  Next, consider a map $\alpha $ of the annulus $S^1\times [5\Phi _-,5\Phi _+]$
into $R^2$, where $\Phi _-$ and $\Phi _+$ are the lower and the upper
bound of the fuction $\Phi $ correspondingly, such that

   $(i)$ The pull-back under $\alpha $ of the form
$dx\wedge dy$ on $R^2$ equals $-d\Phi \wedge d\Theta $.

   $(ii)$ The map $\alpha $ is bijective on $I\times [5\Phi _-,5\Phi _+]$
where $I\subset S^1$ is some closed subset,
such that $I\cup I_0\cup I_1=S^1$; furthermore, the origin
$0\in R^2$ is a unique double point of the map $\alpha $ on
$S^1\times 0$, that is
$$0=\alpha (0,0)=\alpha (\pm 1,0),$$
and
$\alpha $ is injective on $S^1=S^1\times 0$ minus $\{ 0,\pm 1\}$.

   $(iii)$ The curve $S^1_0=\alpha (S^1\times 0)\subset R^2$ ``bounds''
zero area in $R^2$, that is $\int _{S^1_0}xdy=0$, for the $1-$form
$xdy$ on $R^2$.
\begin{Proposition}
Let $V'$, $W'$ and $F'$ as above. Then there exists
an exact Lagrangian embedding $F:W'\times S^1\to V'\times R^2$
given by $F(w',\Theta )=(F'(w',\rho (\Theta )),\alpha (\Theta ,\Phi ))$.
\end{Proposition}
Proof. We follow as in \cite[2.3$B_3'$]{gro}. Now let
$F^*:W'\times S^1\to V'\times R^2$ be given by
$(w',\Theta )\to (F'(w,\rho (\Theta )),\alpha (\Theta ,\Phi ))$.
Then

   $(i)'$ The pull-back under $F^*$ of the form
$\omega =\omega '+dx\wedge dy$ equals
$d\tilde {l}^*-d\Phi \wedge d\Theta =0$ on $W'\times S^1$.

   $(ii)'$ The set of double points of $F^*$ is
$W'_0\cap W_1'\subset V'=V'\times 0\subset V'\times R^2$.

   $(iii)'$  If $F^*$ has no double point then the
Lagrangian submanifold $W=F^*(W'\times S^1)\subset
(V'\times R^2,\omega '+dx\wedge dy)$ is exact if and only if
$W_0'\subset V'$ is such.

   This completes the proof of Proposition 3.1.

\subsection{Formulation of Hilbert bundles}

Let $(\Sigma ,\lambda )$ be a closed $(2n-1)-$ dimensional manifold
with a contact form $\lambda $.
Let $S\Sigma =R\times \Sigma $ and put
$\xi =\ker (\lambda )$.
Let $J'_\lambda $ be
an almost complex structure
on $S\Sigma $
tamed by
the symplectic form $d(e^a\lambda )$.

We define a metric $g_\lambda $ on
$S\Sigma =R\times \Sigma $ by
\begin{equation}
g_\lambda =d(e^a\lambda )(\cdot ,J_\lambda \cdot )
\end{equation}
which is adapted to $J_\lambda $ and $d(e^a\lambda )$ but not
complete.

\vskip 3pt

In the following we denote by $(V',\omega ')=
((R\times \Sigma ),d(e^a\lambda ))$
and
$(V,\omega )=(V'\times R^2,\omega '+dx\wedge dy)$
with the metric $g=g'\oplus g_0 $ induced by
$\omega (\cdot ,J\cdot )$($J=J'\oplus i$ and
$W\subset V$ a Lagrangian submanifold which was constructed in
section 3.1.

   Let $\bar V=D\times V$, then $\pi _1:\bar V\to D$ be
a symplectic vector bundle. Let $\bar J$ be an almost
complex structure on $\bar V$ such that $\pi _1:\bar V\to D$ is
a holomorphic map and each fibre $\bar V_z=\pi _1(z)$ is
a $\bar J$ complex submanifold. Let
$H^k(D)$ be the space of $H^k-$maps
from $D$ to $\bar V$, here
$H^k$ represents Sobolev derivatives up to order $k$. Let
$\bar W=\partial D\times W$, $\bar p=\{1\}\times p$,
$W^{\pm}=\{\pm i\}\times W $ and

$${\cal D}^k
=\{ \bar u \in H^k(D)|
\bar u(x)\in \bar W \ a.e \ for \ x\in \partial D \ and \ \bar u(1)=\bar p,
\bar u(\pm i)\in \{\pm i\}\times W\}$$
for $k\geq 100$.
\begin{Lemma}
Let $W$ be a closed Lagrangian submanifold in
$V$. Then,
${\cal D}^k$
is a pseudo-Hilbert manifold with the tangent bundle
\begin{equation}
T{\cal D}^k
=\bigcup _{\bar u\in {\cal {D}}^k}
\Lambda ^{k-1}
\end{equation}
here
$$\Lambda ^{k-1}=\{ \bar w\in H^{k-1}(\bar u^*(T\bar V)|
\bar w(1)=0, and \ \bar w(\pm i)\in T\bar W\} $$
\end{Lemma}

\begin{Note}
Since $W$ is not regular we know that ${\cal D}^k$ is
in general complete, however it is enough for our purpose.
\end{Note}
Proof: See \cite{al,kl}.

\vskip 3pt

   Now we consider  a section
from ${\cal D}^k$ to
$T{\cal D}^k$ follows as in
\cite{al,gro}, i.e.,
let $\bar \partial :{\cal D}^k\rightarrow T{\cal D}^k$
be the Cauchy-Riemmann section
\begin{equation}
\bar \partial \bar u={{\partial \bar u}\over {\partial s}}
+J{{\partial \bar u}\over {\partial t}}  \label{eq:CR}
\end{equation}
for $\bar u\in {\cal D}^k$.

\begin{Theorem}
The Cauchy-Riemann section $\bar \partial $ defined in (\ref{eq:CR})
is a Fredholm section of Index zero.
\end{Theorem}
Proof. According to the definition of the Fredholm section,
we need to prove that
$\bar u\in {\cal D}^k$, the linearization
$D\bar \partial (\bar u)$ of $\bar \partial $ at $\bar u$ is
a linear Fredholm
operator.
Note that
\begin{equation}
D\bar \partial (\bar u)=D{\bar \partial _{[\bar u]}}
\end{equation}
where
\begin{equation}
(D\bar \partial _{[\bar u]})v=\frac{\partial \bar v}{\partial s}
+J\frac{\partial \bar v}{\partial t}+A(\bar u)\bar v
\end{equation}
with
$$\bar v|_{\partial D}\in (\bar u|_{\partial D})^*T\bar W$$
here $A(\bar u)$ is $2n\times 2n$
matrix induced by the torsion of
almost complex structure, see \cite{al,gro} for the computation.

   Observe that the linearization $D\bar \partial (\bar u)$ of
$\bar \partial $ at $\bar u$ is equivalent to the following Lagrangian
boundary value problem
\begin{eqnarray}
&&{{\partial \bar v}\over {\partial s}}+\bar J
{{\partial \bar v}\over {\partial t}}
+A(\bar u)\bar v=\bar f, \ \bar v\in \Lambda ^k(\bar u^*T\bar V)\cr
&&\bar v(t)\in T_{\bar u(t)}W, \ \ t\in {\partial D}  \label{eq:Lin}
\end{eqnarray}
One
can check that (\ref{eq:Lin})
defines a linear Fredholm operator. In fact,
by proposition 2.2 and Lemma 2.1, since the operator $A(\bar u)$ is a compact,
we know that the operator $\bar \partial $ is a nonlinear Fredholm operator
of the index zero.

\begin{Definition}
Let $X$ be a Banach manifold and $P:Y\to X$ the Banach
vector bundle.
A Fredholm section $F:X\rightarrow Y$ is
proper if $F^{-1}(0)$ is a compact set and is called
generic if $F$ intersects the zero section transversally, see \cite{al,gro}.
\end{Definition}
\begin{Definition}
$deg(F,y)=\sharp \{ F^{-1}(0)\} mod2$ is called the Fredholm
degree of a Fredholm section (see\cite{al,gro}).
\end{Definition}
\begin{Theorem}
Assum that $\bar J=i\oplus J$ on $\bar V$ and $i$ is complex structure
on $D$ and $J$ the almost complex structure on
$V$ which is integrable near point $p$. Then the Fredholm section
$F=\bar \partial : {\cal D}^k\rightarrow T{\cal D}^k$
constructed in (\ref{eq:CR}) has degree one, i.e.,
$$deg(F,0)=1$$
\end{Theorem}
Proof: We assume that $\bar u:D\mapsto \bar V$ be a $\bar J-$holomorphic disk
with boundary $\bar u(\partial D)\subset \bar W$ and
by the assumption that $\bar u$ is homotopic to the
map $\bar u_1=(id,\bar p)$. Since almost complex
structure ${\bar J}$ splits and
is tamed by  the symplectic form $\bar \omega $,
by stokes formula,
we conclude the second component $u: D\rightarrow
V$ is a constant map. Because $u(1)=p$, We know that
$F^{-1}(0)=(id,p)$.
Next we show that the linearizatioon $DF_{(id,p)}$ of $F$ at $(id,p)$ is
an isomorphism from $T_{(id,p)}{\cal D}^k$ to $E$.
This is equivalent to solve the equations
\begin{eqnarray}
{\frac {\partial \bar v}{\partial s}}+J{\frac {\partial \bar v}{\partial t}}=f\\
\bar v|_{\partial D}\subset T_{(id ,p)}\bar W
\end{eqnarray}
here $\bar J=i+J(p)$. By Lemma 2.1,
we know that $DF((id,p))$ is an isomorphism.
Therefore $deg(F,0)=1$.

\section{Anti-holomorphic sections}

In this section we construct a Fredholm
section which is not proper as in
\cite{al,gro}.

  Let $(V',\omega ')=(S\Sigma ,d(e^a\lambda ))$ and
$(V,\omega )=(V'\times C, \omega '\oplus \omega _0)$,
$W$ as in section3 and $J=J'\oplus i$, $g=g'\oplus g_0$,
$g_0$ the standard metric on $C$.

   Now let $c\in C$ be a non-zero vector. We consider
$c$ as an anti-holomorphic homomorphism
$c:TD\to TV'\oplus TC$, i.e., $c({{\partial }\over
{\partial \bar z}})
=(0,c\cdot {{\partial }\over {\partial z}}).$
Since the constant section $c$ is not a section of the
Hilbert bundle in section 3 due to $c$ is not
tangent to the Lagrangian submanifold $W$, we must modify it as follows:

\vskip 3pt

  Let $c$ as above, we define
\begin{eqnarray}
c_{\chi ,\delta }(z,v)=\left\{ \begin{array}{ll}
c \ \ \ &\mbox{if\  $|z|\leq 1-2\delta $,}\cr
0 \ \ \ &\mbox{otherwise}
\end{array}
\right.
\end{eqnarray}
Then by using the cut off function $\varphi _h(z)$ and
its convolution with section
$c_{\chi ,\delta }$, we obtain a smooth section
$c_\delta$ satisfying

\begin{eqnarray}
&&c_{\delta }(z,v)=\left\{ \begin{array}{ll}
c \ \ \ &\mbox{if\  $|z|\leq 1-3\delta $,}\cr
0 \ \ \ &\mbox{if\  $|z|\geq 1-\delta $.}
\end{array}
\right.   \cr
&&|c_\delta |\leq |c|
\end{eqnarray}
for $h$ small enough, for the convolution theory see
\cite[ch1,p16-17,Th1.3.1]{hor}.
Then one can easily check
that $\bar c_\delta =(0,0,c_\delta )$
is an anti-holomorphic section tangent to $\bar W$.

\vskip 3pt

Now we put an almost complex structure $\bar J=i\oplus J$ on the
symplectic fibration $D\times V \to D$ such that $\pi _1:D\times
V\to D$ is a holomorphic fibration and $\pi _1^{-1}(z)$ is an
almost complex submanifold. Let $g=\bar \omega (\cdot ,\bar
J\cdot)$ be the metric on $D\times V$.

Now we consider the
equations
\begin{eqnarray}
&&\bar v=(id ,v)=(id, v',f):D\to D\times V'\times C \cr &&\bar
\partial _{J}v=c_\delta \quad or \cr &&\bar \partial _{J'}v'=0,\bar
\partial f=c_\delta \ on \ D\cr &&v|_{\partial
D}:\partial D\to W\ \ \  \label{eq:4.16}
\end{eqnarray}
here $v$ homotopic to constant map $\{ p\}$ relative to $W$. Note
that $W\subset V\times B_{R}(0)$ for $\pi R^2=2\pi R(\varepsilon
)^2$, here $R(\varepsilon )\to 0$ as $\varepsilon \to 0$ and
$\varepsilon $ as in section 3.1.

\begin{Lemma}
Let $\bar v=(id ,v)$ be the solutions of (\ref{eq:4.16}), then one has
the following estimates
\begin{eqnarray}
E({v})=
\{
\int _D(g'({{\partial {v'}}\over {\partial x}},
{J'}{{\partial {v'}}\over {\partial x}})
+g'({{\partial {v'}}\over {\partial y}},
{J'}{{\partial {v'}}\over {\partial y}}) \nonumber \\
+g_0({{\partial {f}}\over {\partial x}},
{i}{{\partial {f}}\over {\partial x}})
+g_0({{\partial {f}}\over {\partial y}},
{i}{{\partial {f}}\over {\partial y}}))d\sigma \}
\leq 4\pi R(\varepsilon )^2.
\end{eqnarray}
\end{Lemma}
Proof: Since $v(z)=(v'(z),f(z))$ satisfy (\ref{eq:4.16})
and $v(z)=(v'(z),f(z))\in V'\times C$
is homotopic to constant map $v_0:D\to \{ p\}\subset W$
in $(V,W)$, by the Stokes formula
\begin{equation}
\int _{D}v^*(\omega '\oplus \omega _0)=0
\end{equation}
Note that the metric $g$ is adapted to the symplectic form
$\omega $ and $J$, i.e.,
\begin{equation}
g=\omega  (\cdot ,J\cdot )
\end{equation}
By the simple algebraic computation, we have
\begin{equation}
\int _{D}{v}^*\omega  ={{1}\over {4}}
\int _{D^2}(|\partial v|^2
-|\bar {\partial }v|^2)=0
\end{equation}
and
\begin{equation}
|\nabla v|={{1}\over {2}}(
|\partial v|^2 +|\bar \partial v|^2
\end{equation}
Then
\begin{eqnarray}
E(v)&=&\int _{D} |\nabla v| \nonumber \\
      &=&\int _{D}\{ {{1}\over {2}}(
|\partial v|^2+|\bar \partial v|^2)\})d\sigma \nonumber \\
&=&\int _D|c_\delta |_{\bar g}^2d\sigma
\end{eqnarray}
By Cauchy integral formula,
\begin{equation}
f(z)={{1}\over {2\pi i}}\int _{\partial D}{{f(\xi)}\over
{\xi -z}}d\xi
+{{1}\over {2\pi i}}
\int _D{{\bar \partial f(\xi )}\over {\xi -z}}d\xi\wedge d\bar \xi
\label{eq:Ch}
\end{equation}
Since $f$ is smooth up to the boundary,
we integrate the two sides on $D_r$ for $r<1$, one get
\begin{eqnarray}
\int _{\partial D_r}f(z)dz&=&
\int _{\partial D_r}{{1}\over {2\pi i}}\int _{\partial D}{{f(\xi)}\over
{\xi -z}}d\xi dz
+\int _{\partial D_r}{{1}\over {2\pi i}}
\int _D{{\bar \partial f(\xi )}\over {\xi -z}}d\xi\wedge d\bar \xi \cr
&=&0+{{1}\over {2\pi i}}
\int _D\int _{\partial D_r}{{\bar \partial
f(\xi )}\over {\xi -z}}dzd\xi\wedge d\bar \xi \cr
&=&{{1}\over {2\pi i}}
\int _D2\pi i\bar \partial
f(\xi )d\xi\wedge d\bar \xi \label{eq:Ch1}
\end{eqnarray}
Let $r\to 1$, we get

\begin{eqnarray}
\int _{\partial D}f(z)dz=
\int _D\bar \partial
f(\xi )d\xi\wedge d\bar \xi \label{eq:Ch2}
\end{eqnarray}

By the equations (\ref{eq:4.16}),
one get
\begin{equation}
\bar \partial f=c \ on \ D_{1-2\delta }
\end{equation}
So, we have
\begin{equation}
2\pi i(1-2\delta )c=\int _{\partial D}f(z)dz-
\int _{D-D_{1-2\delta }}\bar \partial
f(\xi )d\xi\wedge d\bar \xi \label{eq:Ch3}
\end{equation}
So,
\begin{eqnarray}
|c|&\leq &{{1}\over {2\pi (1-2\delta )}}|\int _{\partial D}f(z)dz|
+|\int _{D-D_{1-2\delta }}\bar \partial
f(\xi )d\xi\wedge d\bar \xi |\cr
&\leq &{{1}\over {2\pi (1-2\delta )}}
2\pi |diam(pr_2(W))+c_1c_2|c|(\pi -\pi (1-2\delta )^2))
\label{eq:Ch4}
\end{eqnarray}
Therefore, one has
\begin{eqnarray}
|c|&\leq &
c(\delta )R(\varepsilon )\label{eq:Ch5}
\end{eqnarray}
and
\begin{eqnarray}
E(v)&=&\pi \int _D|c_\delta |_{\bar g}^2  \cr
    &=&\pi c(\delta )^2R(\varepsilon )^2.
\end{eqnarray}
This finishes the proof of Lemma.

\begin{Proposition}
For $|c|\geq 2c(\delta )R(\varepsilon )$, then the
equations (\ref{eq:4.16})
has no solutions.
\end{Proposition}
Proof. By \ref{eq:Ch5}, it is obvious.

\begin{Theorem}
The Fredholm section $F_1=\bar \partial _{\bar J}+\bar c_\delta
: {\cal  {D}}^k\rightarrow E$ is not proper.
\end{Theorem}
Proof. By the Proposition 4.1 and Theorem 3.2,
it is obvious(see\cite{al, gro}).

\section{$J-$holomorphic section}

Recall that $W(-K,K)\subset W\subset V'\times R^{2}$ as in section
3. The Riemann metric $g$ on $V'\times R^{2}$ induces a metric
$g|W$.

   Now let $c\in C$ be a non-zero vector and
$c_\delta $ the induced anti-holomorphic section. We consider the
nonlinear inhomogeneous equations (\ref{eq:4.16}) and transform it
into $\bar J-$holomorphic map by considering its graph as in
\cite[p319,1.4.C]{gro} or \cite[p312,Lemma5.2.3]{al}.

Denote by $Y^{(1)}\to D\times V$ the bundle of homomorphisms
$T_s(D)\to T_v(V)$. If $D$ and $V$ are given the disk and the
almost K\"ahler manifold, then we distinguish the subbundle
$X^{(1)}\subset Y^{(1)}$ which consists of complex linear
homomorphisms and we denote $\bar X^{(1)}\to D\times V$ the
quotient bundle $Y^{(1)}/X^{(1)}$. Now, we assign to each
$C^1$-map $ v:D\to V$ the section $\bar \partial v$ of the bundle
$\bar X^{(1)}$ over the graph $\Gamma _v\subset D\times V$ by
composing the differential of $v$ with the quotient homomorphism
$Y^{(1)}\to \bar {X}^{(1)}$. If $c_\delta :D\times V\to \bar X$ is
a $H^k-$ section we write $\bar \partial v=c_\delta $ for the
equation $\bar \partial v=c_\delta |\Gamma _v$.

\begin{Lemma}
(Gromov\cite[1.4.$C'$]{gro})There exists a unique almost complex
structure $J_g$ on $D\times V$(which also depends on the given
structures in $D$ and in $V$), such that the (germs of)
$J_\delta-$holomorphic sections $v:D\to D\times V$ are exactly and
only the solutions of the equations $\bar \partial v=c_\delta $.
Furthermore, the fibres $z\times V\subset D\times V$ are
$J_\delta-$holomorphic( i.e. the subbundles $T(z\times V)\subset
T(D\times V)$ are $J_\delta-$complex) and the structure
$J_\delta|z\times V$ equals the original structure on $V=z\times
V$. Moreover $J_\delta $ is tamed by $k\omega _0\oplus \omega $
for $k$ large enough which is independent of $\delta $.
\end{Lemma}

\section{Gromov's $C^0-$convergence theorem}

\subsection{Analysis of Gromov's figure eight}

Since $W'\subset S\Sigma $ is an exact Lagrangian submanifold and
$F'_t$ is an exact Lagrangian isotopy(see section 3.1). Now we
carefully check the Gromov's construction of Lagrangian
submanifold $W\subset V'\times R^2$ from the exact Lagrangian
isotopy of $W'$ in section 3.

   Let $S^1\subset T^*S^1$ be a zero section and $S^1=\cup _{i=1}^4S_i$
be a partition of the zero section $S^1$ such that $S_1=I_0$,
$S_3=I_1$. Write $S^1\setminus \{I_0\cup I_1\}=I_2\cup I_3$ and
$I_0=(-\delta ,-{{5}\over {6}}\delta ]\cup (-{{5\delta }\over
{6}},+{{5\delta }\over {6}})\cup [{{5\delta }\over {6}},
\delta)=I_0^-\cup I_0'\cup I_0^+$, similarly $I_1=(1-\delta
,1-{{5}\over {6}}\delta ]\cup (1-{{5\delta }\over {6}},1+{{5\delta
}\over {6}})\cup [1+{{5\delta }\over {6}}, 1+\delta)=I_1^-\cup
I_1'\cup I_1^+$. Let $S_2=I_0^+\cup I_2\cup I_1^-$, $S_4=I_1^+\cup
I_3\cup I_0^+$. Moreover, we can assume that the double points of
map $\alpha $ in Gromov's figure eight is contained in $(\bar
I_0'\cup \bar I_1')\times [\Phi _-,\Phi _+]$, here $\bar
I_0'=(-{{5\delta }\over {12}},+{{5\delta }\over {12}})$ and $\bar
I_1'=(1-{{5\delta }\over {12}},1+{{5\delta }\over {12}})$. Recall
that $\alpha : (S^1\times [5\Phi _-,5\Phi _+])\to R^2$ is an exact
symplectic immersion, i.e., $\alpha ^*(-ydx)-\Psi d\Theta =dh$,
$h:T^*S^1\to R$. By the construction of figure eight, we can
assume that $\alpha '_i=\alpha |((S^1\setminus I'_i)\times [5\Phi
_-,5\Phi _+])$ is an embedding for $i=0,1$.
 Let $Y=\alpha (S^1\times
[5\Phi _-,5\Phi _+])\subset R^2$ and $Y_i=\alpha (S_i\times [5\Phi
_-,5\Phi _+])\subset R^2$. Let $\alpha _i=\alpha |Y_i(S^1\times
[5\Phi _-,5\Phi _+])$. So, $\alpha _{i}$ puts the function $h$ to
the function $h_{i0}=\alpha _{i}^{-1*}h$ on $Y_i$. We extend the
function $h_{i0}$ to whole plane $R^2$. In the following we take
the liouville form $\beta _{i0}=-ydx -dh_{i0} $ on $R^2$. This
does not change the symplectic form $dx\wedge dy $ on $R^2$. But
we have $\alpha _i^*\beta =\Phi d\Theta $ on $(S_i\times [5\Phi
_-,5\Phi _+])$ for $i=1,2,3,4$.  Finally, note that
\begin{eqnarray}
&&F:W'\times S^1\to V'\times R^2;\cr &&F(w',\Theta )=(F'_{\rho
(\Theta )}(w'),\alpha (\Theta ,\Phi (w',\rho (\Theta )).
\end{eqnarray}
Since $\rho (\Theta )=0 $ for $\Theta \in I_0$ and $\rho (\Theta
)=1$ for $\Theta \in I_1$, we know that $\Phi (w',\rho (\Theta
))=0$ for $\Theta \in I_0\cup I_1$. Therefore,
\begin{eqnarray}
F(W'\times I_0)=W'\times \alpha (I_0); F(W'\times I_1)=W'\times
\alpha (I_1).
\end{eqnarray}

\subsection{Gromov's Schwartz lemma}

In our proof
we need a crucial
tools, i.e., Gromov's Schwartz Lemma as in \cite{gro}.
We first consider the case without boundary.
\begin{Proposition}
Let $(V,J,\mu )$ be as in section 4 and
$V_K$ the compact part of $V$.  There exist constants
$\varepsilon _0$ and $C$(depending only on the $C^0-$ norm of $\mu $ and on the
$C^\alpha $ norm of $J$ and $A_0$) such
that every $J-$holomorphic map of the unit disc to an $\varepsilon _0$-ball
of $V$ with center in $V_K$ and area less than $A_0$ has its derivatives
up to order $k+1+\alpha $ on $D_{{1}\over {2}}(0)$
bounded by $C$.
\end{Proposition}
For a proof, see\cite{gro}.

\vskip 3pt

Now we consider the Gromov's Schwartz Lemma
for $J-$holomorphic map with boundary in
a closed Lagrangian submanifold as in
\cite{gro}.

\begin{Proposition}
Let $(V,J,\mu )$ as above and $L\subset V$ be a closed
Lagrangian submanifold and $V_K$ one
compact part of $V$. There exist constants
$\varepsilon _0$ and $C$(depending only on the $C^0-$ norm of $\mu $ and on the
$C^\alpha $ norm of $J$ and $K,A_0$) such
that every $J-$holomorphic map of the half unit disc $D^+$ to a
$\varepsilon _0$-ball
of $V$ with boundary in $L$ and area less than $A_0$ has its derivatives
up to order $k+1+\alpha $ on $D_{{1}\over {2}}^+(0)$
bounded by $C$.
\end{Proposition}
For a proof see \cite{gro}.

\vskip 3pt

Since in our case $W$ is a non-compact Lagrangian submanifold,
Proposition 6.2 can not be used directly but the proofs of
Proposition 6.1-2 still holds in our case.

\begin{Lemma}
Recall that $V=V'\times R^2$. Let $(V,J,\mu )$ as above and
$W\subset V$ be as above and $V_c$ the compact set in $V$. Let
$\bar V=D\times V$, $\bar W=\partial D\times W$, and $\bar
V_c=D\times V_c$. Let $Y=\alpha (S^1\times [5\Phi _-,5\Phi
_+])\subset R^2$. Let $Y_i=\alpha (S_i\times [5\Phi _-,5\Phi
_+])\subset R^2$. Let $\{ X_j\}_{j=1}^q$ be a Darboux covering of
$\Sigma $ and $V'_j=R\times X_j$.
 Let $\partial D=S^{1+}\cup S^{1-}$. There exist
constant $c_0$ such that every $J-$holomorphic map $v$ of the half
unit disc $D^+$ to the $D\times V_j'\times R^2$ with its boundary
$v((-1,1))\subset (S^{1\pm})\times F({\cal L}\times R\times
S_i)\subset \bar W,i=1,..,4$ has
\begin{equation}
area(v(D^+))\leq c_0l^2(v(\partial 'D^+)).
\end{equation}
here $\partial 'D^+=\partial D\setminus [-1,1]$ and $l(v(\partial
'D^+))=length(v(\partial 'D^+))$.
\end{Lemma}
Proof. Let $\bar W_{i\pm }=S^{1\pm}\times F(W'\times S_i)$. Let
$v=(v_1,v_2):D^+\to \bar V=D\times V$ be the $J-$holomorphic map
with $v(\partial D^+)\subset \bar W_{i\pm}\subset \partial D\times
W$, then

\begin{eqnarray}
area(v)&=&\int _{D^+}v^*d(\alpha _0\oplus \alpha )\cr
&=&\int_{D^+}dv^*(\alpha _0\oplus \alpha )\cr &=&\int _{\partial
D^+}v^*(\alpha _0 \oplus \alpha )\cr &=&\int _{\partial
D^+}v_1^*\alpha _0+\int _{\partial D^+}v_2^*\alpha \cr &=&\int
_{\partial 'D^+\cup [-1,+1]}v_1^*\alpha _0 +\int _{\partial
'D^+\cup [-1,+1]}v_2^*(e^a\lambda-ydx-dh_{i0})\cr &=&\int
_{\partial 'D^+\cup [-1,+1]}v_1^*\alpha _0 + \int _{\partial
'D^+}v_2^*(e^a\lambda-ydx-dh_{i0})+B_1 ,\label{eq:al1}
\end{eqnarray}
here $B_1=\int _{[-1,+1]}{v_2}^*(-d\Psi ')$. Now take a zig-zag
curve $C$ in $V_j'\times Y_i$ connecting $v_2(-1)$ and $v_2(+1)$
such that

\begin{eqnarray}
\int _{C} (e^a\lambda+ydx)&=&B_1\cr length(C)&\leq&
k_1length(v_2(\partial 'D^+)) \label{eq:m_1}
\end{eqnarray}
Now take a minimal surface $M$ in $V'_j\times R^2$ bounded by
$v_2(\partial 'D^+)\cup C$, then by the isoperimetric
ineqality(see[\cite[p283]{grob}), we get

\begin{eqnarray}
area(M)&\leq &m_1length (C+v_2(\partial 'D^+))^2\cr &\leq &
m_2length (v_2(\partial 'D^+))^2,
\end{eqnarray}
here we use the (\ref{eq:m_1}).

Since $area(M)\geq \int _M\omega $ and $\int _M\omega =\int
_{D^+}v_2^*\omega =area(v)$, this proves the lemma.

\begin{Lemma}
Let $v$ as in Lemma 6.1, then we have
\begin{equation}
area(v(D^+)\geq c_0(dist(v(0),v(\partial 'D^+)))^2,
\end{equation}
here $c_0$ depends only on $\Sigma ,J,\omega,...,$etc, not on $v$.
\end{Lemma}
Proof. By the standard argument as in \cite[p79]{al}.

\vskip 3pt

The following estimates is a crucial step in our proof.

\begin{Lemma}
Recall that $V=V'\times R^2$. Let $(V,J,\mu )$ as above and
$W\subset V$ be as above and $V_c$ the compact set in $V$. Let
$\bar V=D\times V$, $\bar W=\partial D\times W$, and $\bar
V_c=D\times V_c$. Let $Y=\alpha (S^1\times [5\Phi _-,5\Phi
_+])\subset R^2$. Let $Y_i=\alpha (S_i\times [5\Phi _-,5\Phi
_+])\subset R^2$. Let $\partial D=S^{1+}\cup S^{1-}$. There exist
constant $c_0$ such that every $J-$holomorphic map $v$ of the half
unit disc $D^+$ to the $D\times V'\times R^2$ with its boundary
$v((-1,1))\subset (S^{1\pm})\times F({\cal L}\times R\times
S_i)\subset \bar W,i=1,..,4$ has
\begin{equation}
area(v(D^+))\leq c_0l^2(v(\partial 'D^+)).
\end{equation}
here $\partial 'D^+=\partial D\setminus [-1,1]$ and $l(v(\partial
'D^+))=length(v(\partial 'D^+))$.
\end{Lemma}
Proof. We first assume that $\varepsilon $ in section 3.1 is small
enough. Let $l_0$ is a constant small enough. If $length(\partial
'D^+)\geq l_0$, then Lemma 6.3 holds. If $length(\partial
'D^+)\leq l_0$ and $v(D^+)\subset D\times V_j'\times R^2$, then
Lemma6.3 reduces to Lemma6.1. If $length(\partial 'D^+)\leq l_0$
and $v(D^+)\bar \subset D\times V_j'\times R^2$, then Lemma6.2
imples $area(v)\geq \tau _0>100\pi R(\varepsilon )^2$, this is a
contradiction. Therefore we  proved the lemma.

\begin{Proposition}
Let $(V,J,\mu )$ and $W\subset V$ be as in section 4 and $V_K$ the
compact part of $V$. Let $\bar V$, $\bar V_K$ and $\bar W$ as
section 5.1. There exist constants $\varepsilon _0$ (depending
only on the $C^0-$ norm of $\mu $ and on the $C^\alpha $ norm of
$J$) and $C$(depending only on the $C^0$ norm of $\mu $ and on the
$C^{k+\alpha }$ norm of $J$) such that every $J-$holomorphic map
of the half unit disc $D^+$ to the $D\times V'\times R^2$ with its
boundary $v((-1,1))\subset (S^{1\pm})\times F({\cal L}\times
R\times S_i)\subset \bar W,i=1,..,4$ has its derivatives up to
order $k+1+\alpha $ on $D_{{1}\over {2}}^+(0)$ bounded by $C$.
\end{Proposition}
Proof. One uses Lemma 6.3 and Gromov's proof on Schwartz lemma to
yield proposition 6.3.

\subsection{Removal singularity of $J-$curves}

In our proof we need another crucial tools, i.e., Gromov's removal
singularity theorem\cite{gro}. We first consider the case without
boundary.
\begin{Proposition}
Let $(V,J,\mu )$ be as in section 4 and $V_K$ the compact part of
$V$. If $v:D\setminus \{0\}\to V_K$ be a $J-$holomorphic disk with
bounded energy and bounded image, then $v$ extends to a
$J-$holomorphic map from the unit disc $D$ to $V_K$.
\end{Proposition}
For a proof, see\cite{gro}.

\vskip 3pt

Now we consider the Gromov's removal singularity theorem for
$J-$holomorphic map with boundary in a closed Lagrangian
submanifold as in \cite{gro}.

\begin{Proposition}
Let $(V,J,\mu )$ as above and $L\subset V$ be a closed Lagrangian
submanifold and $V_K$ one compact part of $V$. If $v:(D^+\setminus
\{0\},\partial ''D^+\setminus \{0\})\to (V_K,L)$ be a
$J-$holomorphic half-disk with bounded energy and bounded image,
then $v$ extends to a $J-$holomorphic map from the half unit disc
$(D^+,\partial ''D^+)$ to $(V_K,L)$.
\end{Proposition}
For a proof see \cite{gro}.

\vskip 3pt

\begin{Proposition}
Let $(V,J,\mu )$ and $W\subset V$ be as in section 4 and $V_c$ the
compact set in $V$. Let $\bar V=D\times V$, $\bar W=\partial
D\times W$, and $\bar V_c=D\times V_c$. Then every $J-$holomorphic
map $v$ of the half unit disc $D^+\setminus \{0\}$ to the $\bar V$
with center in $\bar V_c$ and its boundary $v((-1,1)\setminus
\{0\})\subset (S^{1\pm})\times F({\cal L}\times [-K,K]\times
S_i)\subset \bar W$ and
\begin{equation}
area(v(D^+\setminus \{0\}))\leq E
\end{equation}
extends to a $J-$holomorphic map $\tilde v:(D^+, \partial ''D)\to
(\bar V_c, \bar W)$.
\end{Proposition}
Proof. This is ordinary Gromov's removal singularity theorem by
$K-$assumption.

\subsection{$C^0-$Convergence Theorem}

We now recall that the well-known Gromov's compactness
theorem for cusp's curves for the compact symplectic manifolds with
closed Lagrangian submanifolds in it.
For reader's convenience, we first recall the ``weak-convergence''
for closed curves.

\vskip 3pt

{\bf Cusp-curves.}
Take a system of disjoint simple closed curves $\gamma _i$
in a closed surface $S$ for $i=1,...,k$, and denote by $S^0$ the
surface obtained from $S\setminus \cup _{i=1}^k\gamma _i$. Denote by
$\bar S$ the space obtained from $S$ by shrinking every $\gamma _i$
to a single point and observe the obvious map $\alpha :S^0\to \bar S$ gluing
pairs of points $s'_i$ and $s''_i$ in $S^0$, such that
$\bar {s}_i=\alpha (s_i')=\alpha (s_i'')\in \bar S$
are singular (or cuspidal) points in $\bar S$(see\cite{gro}).

\smallskip

An almost complex structure in $\bar {S}$ by definition is that in $S^0$.
A continuous map $\beta :\bar {S}\to V$ is called a (parametrized
$J-$holomorphic) cusp-curve in $V$ if
the composed map $\beta \circ \alpha :S^0\to V$ is
holomorphic.

\vskip 3pt

{\bf Weak convergence.} A sequence of closed $J-$curves
$C_j\subset V$ is said to weakly converge to a cusp-curve
$\bar {C}\subset V$ if the following four conditions are satisfied

(i) all curves $C_j$ are parametrized by a fixed surface $S$
whose almost complex structure depends on $j$, say
$C_j=f_j(S)$ for some holomorphic maps
$$f_j:(S,J_j)\to (V,J)$$

(ii) There are disjoint simple closed curves $\gamma _i\in S$,
$i=1,...,k$, such that $\bar {C}=\bar {f}(\bar {S})$ for
a map $\bar {f}:\bar {S}\to V$ which
is holomorphic for some almost complex
structure $\bar {J}$ on $\bar {S}$.

(iii) The structures $J_j$ uniformly $C^\infty -$converge
to $\bar J$ on compact
subsets in $S\setminus \cup _{i=1}^k\gamma _i$.

(iv) The maps $f_j$ uniformly $C^\infty -$converge to
$\bar f$ on compact subsets in $S\setminus \cup _{i=1}^k\gamma _i$.
Moreover, $f_j$ uniformly $C^0-$converge on entire $S$
to the composed map $S\to
\bar {S}\stackrel {\bar {f}}{\rightarrow }V$.
Furthermore,

$$Area_{\mu }f_j(S)\to Area _{\mu }\bar {f}(\bar {S})\ for \
j\to \infty ,$$
where $\mu $ is a Riemannian metric in $V$ and where
the area is counted with the geometric multiplicity(see\cite{gro}).

\vskip 3pt

{\bf Gromov's Compactness theorem for closed curves.} Let $C_j$
be a sequence of closed $J-$curves of a fixed
genus in a compact manifold $(V, J, \mu ).$ If the areas of $C_j$
are uniformly bounded,
$$Area _{\mu }\leq A,\ j=1,..,$$
then some subsequence weakly converges to
a cusp-curve $\bar C$ in $V$.

\vskip 3pt

{\bf Cusp-curves with boundary.} Let $T$ be a compact
complex manifold with boundary of dimension $1$(i.e., it has an atlas of
holomorphic charts onto open subsets of $C$ or of a closed half plane).
Its double is a compact Riemann surface $S$ with a natureal anti-holomorphic
involution $\tau $ which exchanges $T$ and $S\setminus T$ while
fixing the boundary $\partial T$. IF$f:T\to V$ is a continous
map, holomorphic in the interior
of $T$, it is convenient to extend $f$ to $S$ by
$$f=f\circ \tau$$
Take a totally real submanifold $W\subset (V,J)$ and consider compact
holomorphic curves $C\subset V$ with boundaries, $(\bar C, \partial \bar
C)\subset (V,W)$, which are,
topologically speaking,
obtained by shrinking to points some (short) closed loops in $C$
and also some (short) segments in $C$ between boundary points.
This is seen by looking
on the double $C\cup _{\partial C}C$.

\vskip 3pt

{\bf Gromov's Compactness theorem for curves with boundary.}
Let $V$ be a closed Riemannian manifold, $W$ a totally real closed
submanifold of $V$. Let $C_j$
be a sequence of $J-$curves with boundary in $W$ of a fixed
genus in a compact manifold $(V, J, \mu )$. If the areas of $C_j$
are uniformly bounded,
$$Area _{\mu }\leq A,\ j=1,..,$$
then some subsequence weakly converges to
a cusp-curve $\bar C$ in $V$.

\vskip 3pt

The proofs of Gromov's compactness theorem can found in
\cite{al,gro}. In our case the Lagrangian submanifold $W$ is not
compact, Gromov's compactness theorem can not be applied directly
but its proof is still effective since the $W$ has the special
geometry. In the following we modify Gromov's proof to prove the
$C^0-$compactness theorem in our case.

\vskip 3pt

  Now we state the $C^0-$convergence theorem in our case.

\begin{Theorem}
Let $(V, J, \omega , \mu )$ and $W$
as in section4.
Let $C_j$
be a sequence of $\bar J_\delta -$holomorphic section
$v_j=(id, ((a_j,u_j),f_j)):D\to D\times V$ with
$v_j:\partial D \to \partial D\times W$ and $v_j(1)=(1,p)\in
\partial D\times W$ constructed from section 4.  Then the areas of $C_j$
are uniformly bounded,i.e.,
$$Area _{\mu }(C_j)\leq A,\ j=1,..,$$
and some subsequence weakly converges to a cusp-section $\bar C$
in $V$(see\cite{al,gro}).
\end{Theorem}
Proof. We follow the proofs in \cite{gro}. Write
$v_j=(id,(a_{j},u_{j}),f_j))$ then $|a_{ij}|\leq a_0$
by the ordinary Monotone inequality of minimal
surface without boundary, see following Proposition 7.1.
Similarly $|f_j|\leq R_1$
by using the fact
$f_j(\partial D)$ is bounded in $B_{R_1}(0)$
and $\int _D|\nabla f_j|\leq 4\pi R^2$ via monotone inequality
for minimal surfaces. So, we assume
that $v_j(D)\subset V_c$ for a compact set $V_c$.

\vskip 2pt

1. {\it Removal of a net}.

\smallskip

1a.  Let $\bar V=D\times V$ and
$v_j$ be regular curves. First we study induced metrics
$\mu _j$ in $v_j$.
We apply the ordinary monotone inequality for minimal surfaces
without boundary to small concentric balls
$B_\varepsilon \subset (A_j,\mu _j)$ for
$0<\varepsilon \leq \varepsilon _0$ and conclude by
the standard argument to the inequality
$$Area (B_\varepsilon )\geq \varepsilon ^2, \ for \ \varepsilon
\leq \varepsilon _0;$$
Using this we easily find a interior  $\varepsilon -$net
$F_j\subset (v_j,\mu _j)$ containing $N$ points for
a fixed integer
$N=(\bar V, \bar J,\mu )$, such that every
topological annulus $A\subset v_j\setminus F_j$ satisfies
\begin{equation}
Diam_\mu A\leq 10length_\mu \partial A.
\end{equation}
Furthermore, let $A$ be conformally equivalent to the
cylinder $S^1\times [0,l]$ where $S^1$ is the circle
of the unit length, and let $S_t^1\subset A$ be the curve in $A$ corresponding to
the circle $S^1\times t$for $t\in [0,l]$. Then
obviously
\begin{equation}
\int _a^b(lengthS_t^1)^2dt\leq Area(A)\leq C_5.
\end{equation}
for all $[a,b]\subset [0,l]$. Hence, the annulus $A_t\subset
A$ between the
curves $S^1_t$ and $S^1_{l-t}$ satisfies
\begin{equation}
diam_\mu A_t\leq 20 ({{C_5}\over {t}})
\end{equation}
for all $t\in [0,l]$.

\vskip 3pt

$1b$.  We consider the sets $\partial v_j\cap ((S^{1\pm})\times
F(W'\times I_i^\pm)), i=0,1$. By the construction of Gromov's
figure eight, there exists a finite components, denote it by
\begin{equation}
\partial v_j\cap ((S^{1\pm})\times F({\cal L}\times R\times I_i^\pm))
=\{\bar \gamma _{ij}^{k}\},i=0,1.
\end{equation}
Let $m_i^\pm$ be the middle point of $I_i^\pm$. If
\begin{equation}
\bar \gamma _{ij}^k\cap ((S^{1\pm})\times F({\cal L}\times R\times
m_i^\pm))\ne \emptyset ,i=0,1,
\end{equation}
we choose one point in $\bar \gamma _{ij}^{k}$ as a boundary
puncture point in $\partial v_j$. 
Consider the concentric
$\varepsilon $ half-disks or quadrature $B_\varepsilon (p)$ with
center $p$ on ${\bar \gamma }_{ij}^k$, then 
\begin{equation}
Area(B_\varepsilon (p))\geq \tau _0.
\end{equation}
Since $Area(v_j)\leq E_0$, there exists a uniform finite puncture
points.

Consider the concentric
$\varepsilon $ half-disks or quadrature $B_\varepsilon (p)$ with
center $p$ on $\partial v_j$ and 
\begin{equation}
Area(B_\varepsilon (p))\geq \tau _0,
\end{equation}
we puncture one point on such half-disk or quadrature. 
Since $Area(v_j)\leq E_0$, there exists a uniform finite puncture
points.

So, we find a boundary net $G_j\subset \partial v_j$ containing
$N_1$ points for a fixed integer $N_1(\bar V, \bar J,\mu )$, such
that every topological quadrature or half annulus $B\subset
v_j\setminus \{F_j,G_j\}$ satisfies
\begin{eqnarray}
&&\partial ''B=\partial B\cap \bar W \subset (S^{1\pm})\times
F({\cal L}\times R\times S_i), i=1,2,3,4.
\end{eqnarray}

\vskip 2pt

2. {\it Poincare's metrics}. 2a. Now, let $\mu _j^*$ be a metric
of constant curvature $-1$ in $v_j(D)\setminus F_j\cup G_j$
conformally equivalent to $\mu _j$. Then for every $\mu _j^*-$ball
$B_\rho $ in $v_j\setminus F_j\cup G_j$ of radius $\rho \leq 0.1$,
there exists an annulus $A$ contained in $v_j\setminus F_j\cup
G_j$ such that $B_\rho \subset A_t$ for $t=0.01|log|$(see Lemma
3.2.2in \cite[chVIII]{al}). This implies with $(6.3)$ the uniform
continuity of the (inclusion) maps $(v_j\setminus F_j,\mu _j^*)\to
(\bar V,\bar \mu )$, and hence a uniform bound on the $r^{th}$
order differentials for every $r=0,1,2,...$.

2b. Similarly,
for every $\mu _j^*-$half ball $B_\rho ^+$ in
$v_j\setminus F_j\cup G_j$ of radius
$\rho \leq 0.1$, there
exists a half annulus or quadrature $B$ contained in $v_j\setminus F_j\cup G_j$ such that
$B_\rho ^+\subset B$ with
\begin{eqnarray}
&&\partial ''B=\partial B\cap \bar W \subset (S^{1\pm})\times
F({\cal L}\times R\times S_i), i=1,2,3,4.
\end{eqnarray}
Then, by Gromov's Schwartz Lemma, i.e., Proposition 6.1-6.3
implies the uniform bound on the $r^{th}$ order differentials for
every $r=0,1,2,...$.

\vskip 2pt

3. {\it Convergence of metrics}. Next, by the
standard (and obvious ) properties of hyperbolic surfaces there is a
subsequence(see\cite{al}), which is still denoted by $v_j$, such that

\vskip 3pt

$(a)$. There exist $k$ closed geodesics or geodesic arcs with boundaries in
$\partial v_j\setminus F_j$, say
$$\gamma _i^j\subset (v_j\setminus F_j,\mu ^*_j), i=1,...,k,j=1,2,...,$$
whose $\mu _j^*-$length converges to zero as $j\to \infty $, where
$k$ is a fixed integer.

\vskip 2pt

$(b)$. There exist $k$ closed curves or geodesic arcs
with boundaries in $\partial S$ of  a fixed surface, say $\gamma _j$ in $S$,
and an almost complex structure $\bar J$ on the corresponding
(singular)  surface $\bar S$, such that
the almost complex structure $J_j$ on $v_j\setminus F_j$ induced from
$(V,J)$ $C^\infty -$converge to $\bar J$
outside $\cup _{j=1}^k\gamma _j$. Namely, there exist continuous maps
$g^j:v_j\to \bar S$ which are homeomorphisms
outside the geodesics $\gamma _i^j$, which pinch these
geodesics to the corresponding singular points of
$\bar S$(that are the images of $\gamma _i$) and which
send $F_j$ to a fixed subset
$F$ in the nonsingular locus of $\bar S$. Now,
the convergence $J_j\to \bar J$ is understood as the uniform
$C^\infty -$convergence $g^j_*(J_j)\to \bar J$ on the
compact subsets in the
non-singular  locus $\bar S^*$ of
$\bar S$ which is identified with $S\setminus \cup _{i=1}^k\gamma _i$.

4. {\it $C^0-$interior convergence}.  The limit cusp-curve
$\bar v:\bar S ^*\to  \bar V$, that is a holomorphic map
which is constructed by first taking the maps
$$\bar v_j=(g_j)^{-1}:S\setminus \cup _{i=1}^k\gamma _i\to \bar V$$
Near the nodes of $\bar S$ including
interior nodes and boundary nodes, by the
properties of hyperbolic metric $\mu ^*$ on $\bar S$, the neighbourhoods
of interior nodes are corresponding to the annulis
of the geodesic cycles. By
the reparametrization of $v_j$, called $\bar v_j$ which
is defined on $S$ and extends the maps $\bar v_j:S\to S_j\to V$(see\cite{al,gro}).
Now let $\{ z_i|i=1,...,n\}$ be the interior nodes of $\bar S$. Then
the arguments in \cite{al,gro} yield
the $C^0-$interir convergece near $z_i$.

5. {\it $C^0-$boundary convergence}. Now it is possible that the
boundary of the cusp curve $\bar v$ does not remain in $\bar W$.
Write $\bar v(z)=(h(z),(a(z),u(z)),f(z))$, here $h(z)=z$ or
$h(z)\equiv z_i$, $i=1,...,n$, $z_i$ is cusp-point or bubble
point. We can assume that $\bar p=(1,p)\in \bar v_n$ is a puncture
boundary point. Let $\bar v_1$ be the component of $\bar v$ which
through the point $\bar p$. Let $D=\{z|z=re^{i\theta},0\leq
r\leq1,0\leq \theta \leq 2\pi\}$. We assume that $\bar
v_1:D\setminus \{e^{i\theta _i}\}_{i=1}^k\to V_c$, here
$e^{i\theta _i}$ is node or puncture point. Near $e^{i\theta _i}$,
we take a small disk $D_i$ in $D$ containing only one puncture or
node point $e^{i\theta _i}$. By the reparametrization and the
convergence procedure, we can assume that $\bar v_{1i}=(\bar
v_1|D_i)$ as a map from $D^+\setminus \{0\}\to V_c$ with $\bar
v_1([-1,1]\setminus \{0\})\subset S^1\times F(W'\times S^1)$ and
$area (\bar v_{1i})\leq a_0$, $a_0$ small enough.
 Since $Area(\bar v_{1i})\leq a_0$, there exist curves $c_k$
near $0$ such that $l(\bar v_{1i}(c_k))\leq \delta _1$. By the
construction of convergence, we can assume that $l(\bar
v_n(c_k))\leq 2\delta _1$. If $\bar v_{1i}(\partial c_k)\subset
(S^{1})\times F({\cal L}\times [-N_0,N_0]\times S^1)$, we have
$\bar v_n(\partial c_k)\subset (S^{1})\times F({\cal L}\times
[-2N_0,2N_0]\times S^1)$ for $n$ large enough. Now $\bar v_n(c_k)$
cuts $\bar v_n(D)$ as two parts, one part corresponds to $\bar
v_{1i}$, say $\bar u_n(D)$. Then $area (\bar u_n(D))=
area(h_{n1})+|\Psi '(u_{n2}(c_k^1))-\Psi '(u_{n2}(c_k^2))|$, here
$\partial c_k=\{c_k^1,c_k^2\}$. Then by the proof of Lemma6.1-6.3,
we know that $\bar u_n(\partial D\setminus c_k)\subset
(S^{1})\times F({\cal L}\times [-100N_0,100N_0]\times S^1)$. So,
$\bar v_{1i}([-1,1]\setminus \{0\})\subset S^1\times F({\cal
L}\times [-100N_0,100N_0]\times S^1)$. By proposition 6.6, one
singularity of $\bar v_1$ is deleted. We repeat this procedure, we
proved that $\bar v_1$ is extended to whole $D$. So, the boundary
node or puncture points of $\bar v$ are removed. Then by choosing
the sub-sub-sequences of $\mu ^*_j$ and $\bar v_j$, we know that
$\bar v_j$ converges to $\bar v$ in $C^0$ near the boundary node
or puncture point. This proved the $C^0-$boundary convergence.
Since $\bar v_j(1)=\bar p$, $\bar p\in \bar v(\partial D)$, $\bar
v(\partial D)\subset \bar W$.

\vskip 2pt

6. {\it Convergence of area}. Finally by the
$C^0-$convergence and $area(v_j)=\int _Dv_j^* \bar \omega $,
one easily deduces
$$area(v(S))=\lim _{j\to \infty }(v_j(S_j)).$$

\subsection{Bounded image of $J-$holomorphic curves in $W$}

\begin{Proposition}
Let $v$ be the solutions of equations (4.16), then
$$d_W(p,v(\partial D^2))
=max\{ d_W(p,q)|q\in f(\partial D^2)\} \leq d_0<+\infty$$
\end{Proposition}
Proof. It follows directly from
Gromov's $C^0-$convergence theorem.

\section{Proof of Theorem 1.1}

\begin{Proposition}
If $J-$holomorphic curves $C\subset \bar V$ with
boundary
$$\partial C\subset D^2\times ([0,\varepsilon ]\times \Sigma )
\times R^2$$
and
$$C\cap (D^2\times (\{-3\}\times \Sigma )
\times R^2)\ne \emptyset $$
Then
$$area(C)\geq 2l_0.$$
\end{Proposition}
Proof. It is obvious by monotone inequality argument for
minimal surfaces.

\begin{Note}
we first observe that
any $J-$holomorphic curves with boundary in
$R^+\times \Sigma $ meet the
hypersurface $\{-3\}\times \Sigma $ has
energy at least $2l_0$, so we take
$\varepsilon $ small enough such that
the Gromov's figure eight contained
in $B_{R(\varepsilon )}\subset C$ for $\varepsilon $ small enough and
the energy
of solutions in section 4 is smaller than
$l_0$. we specify the constant
$a_0$, $\varepsilon $ in section 3.1-3 such that
the above conditions satisfied.
\end{Note}

\begin{Theorem}
There exists a non-constant $J-$holomorphic map $u: (D,\partial D)\to
(V'\times C,W)$ with $E(u)\leq 4\pi R(\varepsilon )^2$
for $\varepsilon $ small enough such that $4\pi R(\varepsilon )^2\leq
l_0$.
\end{Theorem}
Proof.  By Proposition 5.1, we know that the image $\bar v(D)$
of solutions of equations (\ref{eq:4.16}) remains a bounded or
compact part of the non-compact Lagrangian submanifold
$W$. Then, all arguments in \cite{al,gro} for the case $W$ is closed
in $S\Sigma\times R^2$ can be extended to our case, especially
Gromov's $C^0-$converngence theorem applies. But the results in
section 4 shows the solutions of equations (\ref{eq:4.16}) must
denegerate to a cusp curves, i.e., we obtain a Sacks-Uhlenbeck-Gromov's
bubble, i.e., $J-$holomorphic sphere or disk with boundary
in $W$, the exactness of $\omega $ rules out
the possibility of $J-$holomorphic sphere. For the more detail, see the proof
of Theorem 2.3.B in \cite{gro}.

\vskip 3pt

{\bf Proof of Theorem 1.1}.
If $(\Sigma ,\lambda )$ has no Reeb chord, then
we can construct a Lagrangian
submanifold $W$ in $V=V'\times C$, see section 3.
Then as in \cite{al,gro}, we construct an
anti-holomorphic section $c$ and
for large vector $c\in C$ we know
that the nonlinear
Fredholm section or Cauchy-Riemann section
has no solution, this implies that the section is non-proper, see section 4.
The non-properness of the
section and
the Gromov's compactness
theorem in section 6 implies the existences of the cusp-curves.
So, we must have the $J-$holomorphic sphere or
$J-$holomorphic disk with bounadry in $W$.
Since the symplectic manifold $V$ is
an exact symplectic mainifold and $W$ is an exact
Lagrangian submanifold in $V$, by Stokes formula,
we know that
the possibility of $J-$holomorphic sphere or disk
elimitated.
So our priori assumption does not holds which
implies the contact maifold $(\Sigma ,\lambda )$
has at least Reeb chord. This finishes
the proof of Theorem 1.1.


\begin{thebibliography}{99}

\bibitem{ab} Abbas, C., Finite energy surface and the chord
probems, Duke. Math.J., Vol.96,No.2,1999,pp.241-316.

\bibitem{ar} Arnold, V. I., First steps in symplectic topology, Russian Math.
Surveys 41(1986),1-21.

\bibitem{ag}  Arnold, V.\& Givental, A., Symplectic Geometry,
in: Dynamical Systems IV, edited by V. I. Arnold and S. P. Novikov,
Springer-Verlag, 1985.

\bibitem{al}  Audin, M\& Lafontaine, J., eds.: Holomorphic Curves in Symplectic
Geometry. Progr. Math. 117, (1994) Birkha\"{u}ser, Boston.

\bibitem{ch}  Chaperon, M., Questions de geometrie symplectique,
in Seminaire Bourbaki, Asterisque 105-106(1983), 231-249.



\bibitem{gi} Givental, A. B., Nonlinear generalization of the
Maslov index, Adv. in Sov. Math., V.1,
AMS, Providence, RI, 1990.

\bibitem{gra} Gray, J.W., Some global properties of contact
structures. Ann. of Math., 2(69): 421-450, 1959.



\bibitem{gro}  Gromov, M., Pseudoholomorphic Curves in Symplectic manifolds.
Inv. Math. 82(1985), 307-347.

\bibitem{grob}  Gromov, M., Partial Differential
Relations, Springer-Verlag, 1986.


\bibitem{ho}  Hofer, H., Pseudoholomorphic curves in symplectizations with
applications to the Weinstein conjection in dimension three.
Inventions Math., 114(1993), 515-563.


\bibitem{hv} Hofer, H.\& Viterbo, C., The Weinstein
conjecture in the presence of holomorphic spheres, Comm. Pure
Appl. Math.45(1992)583-622.


\bibitem{hor} H\"ormander, L., The Analysis of Linear Partial
Differential Operators I, Springer-Verlag, 1983.


\bibitem{kl}  Klingenberg, K., Lectures on closed Geodesics, Grundlehren der Math.
Wissenschaften, vol 230, Spinger-Verlag, 1978.

\bibitem{ls} Lalonde, F \& Sikorav, J.C., Sous-Vari\`t\`es Lagrangiennes
et lagrangiennes exactes des fibr\`es cotangents, Comment. Math. Helvetici 66(1991)
18-33.

\bibitem{ma} Ma, R., Legendrian submanifolds and A Proof on
Chord Conjecture, Boundary Value Problems, Integral Equations and
Related Problems, edited by J K Lu \& G C Wen, World Scientific,
135-142,2000.

\bibitem{mo}  Mohnke, K.: Holomorphic Disks and the Chord Conjecture, Annals of Math., (2001), 154:219-222.


\bibitem{su} Sacks, J. and Uhlenbeck,K., The existence of
minimal 2-spheres. Ann. Math., 113:1-24, 1983.

\bibitem{sm} Smale, S., An infinite dimensional
version of Sard's theorem, Amer. J. Math. 87: 861-866, 1965.


\bibitem{wen}  Wendland, W., Elliptic systems in the plane, Monographs and
studies in Mathematics 3, Pitman, London-San Francisco, 1979.

\end{thebibliography}
\end{document}